\documentclass[12pt,a4paper]{article}
\usepackage[cp1251]{inputenc}
\usepackage[T2A]{fontenc}
\usepackage[english]{babel}

\usepackage{amsmath,amsfonts,amssymb,mathrsfs,amsopn,indentfirst,amscd}

\usepackage{graphicx}
\usepackage{amsthm}
\usepackage{hyperref}

\usepackage{amssymb}
\usepackage{amsmath}

\newcommand{\ver}{{\rm ver}}
\newcommand{\vo}{{\rm vol}}

\newtheorem*{corollary*}{Corollary}

\begin{document}

\title{
On Some Estimate for the Norm \\ of an Interpolation Projector}

\author{Mikhail Nevskii\footnote{ Department of Mathematics,  P.G.~Demidov Yaroslavl State University,\newline Sovetskaya str., 14, Yaroslavl, 150003, Russia, \newline
               mnevsk55@yandex.ru,  orcid.org/0000-0002-6392-7618 } 
               }
      
\date{May 7, 2022}
\maketitle

\begin{abstract}

\smallskip
Let  $Q_n=[0,1]^n$ be the unit cube   in ${\mathbb R}^n$ and let 
$C(Q_n)$ be a space of ~continuous functions 
$f:Q_n\to{\mathbb R}$ with the  norm 
$\|f\|_{C(Q_n)}:=\max_{x\in Q_n}|f(x)|.$ By~$\Pi_1\left({\mathbb R}^n\right)$  denote a set of polynomials of degree
$\leq 1$, i.\,e., a set of linear functions on
${\mathbb R}^n$. The interpolation projector  
$P:C(Q_n)\to \Pi_1({\mathbb R}^n)$ with the~nodes $x^{(j)}\in Q_n$
is defined by the~equalities
$Pf\left(x^{(j)}\right)=
f\left(x^{(j)}\right)$,  $j=1,$ $\ldots,$ $ n+1$. Let $\|P\|_{Q_n}$ be the norm of $P$ as an operator from
$C(Q_n)$ to~$C(Q_n)$. 
If~$n+1$ is~an~Hadamard number, then there exists a nondegenerate  regular simplex having the vertices at vertices of $Q_n$.  We discuss some approaches to~get inequalities 
 of the form $||P||_{Q_n}\leq c\sqrt{n}$ for the norm of the corresponding projector $P$.



\end{abstract}


\section{Introduction}\label{nev_s1}

Let $K$ be a convex body in
${\mathbb R}^n$. i.\,e., a compact convex subset of ${\mathbb R}^n$ with nonempty interior. Denote by   
$C(K)$ a space of continuous functions 
$f:K\to{\mathbb R}$ with  the~uniform norm
$$\|f\|_{C(K)}:=\max\limits_{x\in K}|f(x)|.$$
By $\Pi_1\left({\mathbb R}^n\right)$ we mean a set of polynomials in
$n$ variables of degree
$\leq 1$, i.\,e.,  of  linear functions on ${\mathbb R}^n$.
For $x^{(0)}\in {\mathbb R}^n,  R>0$, by $B(x^{(0)};R)$ we denote the $n$-dimensional Euclidean ball given
by the inequality 
$\|x-x^{(0)}\|\leq R$. Here 
$$\|x\|:=\sqrt{(x,x)}=\left(\sum\limits_{i=1}^n x_i^2\right)^{1/2}.$$
By definition,  $B_n:=B(0;1), Q_n:=[0,1]^n, Q_n^\prime:=[-1,1]^n$. 

Let $S$ be a nondegenerate simplex in  ${\mathbb R}^n$ with the vertices 
$x^{(j)}=\left(x_1^{(j)},\ldots,x_n^{(j)}\right),$ 
$1\leq j\leq n+1.$ Consider the following {\it vertex matrix} of this simplex:
$${\bf S} :=
\left( \begin{array}{cccc}
x_1^{(1)}&\ldots&x_n^{(1)}&1\\
x_1^{(2)}&\ldots&x_n^{(2)}&1\\
\vdots&\vdots&\vdots&\vdots\\
x_1^{(n+1)}&\ldots&x_n^{(n+1)}&1\\
\end{array}
\right).$$
Let ${\bf S}^{-1}$ $=(l_{ij})$. Linear  polynomials 
$\lambda_j(x):=
l_{1j}x_1+\ldots+
l_{nj}x_n+l_{n+1,j}$
have the~property  
$\lambda_j\left(x^{(k)}\right)$ $=$ 
$\delta_j^k$.
We call $\lambda_j$ {\it the basic Lagrange polynomials corresponding to $S$.}
For~an~arbitrary $x\in{\mathbb R}^n$, 
$$x=\sum_{j=1}^{n+1} \lambda_j(x)x^{(j)}, \quad \sum_{j=1}^{n+1} \lambda_j(x)=1.$$
These equalities mean that $\lambda_j(x)$ are {\it  the barycentric coordinates of $x$}. Also we have 
$\lambda_j(x)=\frac{\Delta_j(x)}{\Delta}$,
where $\Delta=\det(\bf S)$ and $\Delta_j(x)$ is obtained from $\Delta$ by replacing the $j$th row by the row $(x_1 \ \ldots \  x_n \ 1)$.
For details, see~\cite[\S 1.1]{nevskii_monograph}.

We say that an interpolation projector
 $P:C(K)\to \Pi_1({\mathbb R}^n)$ corresponds to~a~simplex $S\subset K$  if the nodes of $P$
 coincide with the vertices of $S$. This projector is defined by the~equalities
$Pf\left(x^{(j)}\right)=
f\left(x^{(j)}\right).$
The following analogue of the Lagrange interpolation formula holds:
\begin{equation}\label{interp_Lagrange_formula}
Pf(x)=\sum\limits_{j=1}^{n+1}
f\left(x^{(j)}\right)\lambda_j(x). 
\end{equation}
Denote by $\|P\|_K$ the norm of $P$ as an operator from $C(K)$
in $C(K)$. From (\ref{interp_Lagrange_formula}), it~follows that
\begin{equation}\label{norm_of_P_lambda_j}
\|P\|_K=
\max_{x\in K}\sum_{j=1}^{n+1}
|\lambda_j(x)|=\max_{x\in K}\sum_{j=1}^{n+1}\frac{|\Delta_j(x)|}{|\Delta|}.
\end{equation}
If $K$ is a convex polytope, then the maxima in \eqref{norm_of_P_lambda_j} also can be taken only over $x\in \ver(K)$, where $\ver(K)$ denotes the vertex set of  $K$.

In this paper, we consider the case when $n+1$ is~an~Hadamard number,\linebreak i.\,e., there exists an Hadamard matrix of~order $n+1$. 
{\it An Hadamard matrix of~order~$m$} is~a~square binary matrix $\bf H$ with entries either $1$ or $-1$ which satisfies 
$${\bf H}^{-1}=\frac{1}{m}\, {\bf H}^{T}.$$
This means that the rows of ${\bf H}$ are pairwise orthogonal, with respect to the standard scalar product on ${\mathbb R}^m.$ 

The order of an Hadamard matrix is 1 or 2 or some multiple of 4 (see \cite{hall_1970}).\linebreak  But it is
still unknown whether an Hadamard matrix exists for every order \linebreak  of the form
$m=4k$. This  is one of the longest lasting open
problems in Mathe\-matics
called {\it the Hadamard matrix conjecture}.
The orders below 1500 for which Hadamard matrices are not known are 668, 716, 892, 956, 1132, 1244, 1388, and~1436 (for  links, see, e.g.,~\cite{horadam_2007}, \cite{manjhi_2022}).


Two Hadamard matrices are called equivalent iff one can be obtained  from \linebreak another by 
 performing a finite sequence of the following operations: multiplication of~some rows or columns by $-1$, or permutation of rows or columns. Up~to equi\-valence, there~is a~unique Hadamard matrix of orders 1, 2, 4, 8, and 12. There~are 5 equivalence classes of Hadamard matrices  of order 16, 3 of order 20,  60 of order 24, and 487 of order 28. 
For orders 32, 36, and 40, the number of equivalence classes is~much greater. 
For $ n = 32$ there are at least 3,578,006 equivalence classes; for~$ n = 36$, at~least
4,745,357  (see \cite{horadam_2007}).

Let $Q$ be an arbitrary $n$-dimensional cube. If $n+1$ is an Hadamard number, then there exists a nondegenerate regular simplex having the vertices coinciding with vertices of
 $Q$. In the paper, we discuss some approaches to obtain estimates of~the~form  $||P||_{Q}$ $\leq$  $c\sqrt{n}$ for the norm of the corresponding interpolation projector. 
 For convenience, we will take  $Q=Q_n$  or $Q=Q_n^\prime$.


\section{Making Use of the  Numbers $h_n$}\label{nev_s2}
 Denote by $h_n$ the maximum value of a determinant of order $n$ with entries \linebreak $0$ or~$1$. Let $\nu_n$ be the maximum volume of an $n$-dimensional simplex
 contained in $Q_n$. These~numbers satisfy  the equality $h_n=n!\nu_n$ (see \cite{hudelson_1996}).
For any $n$, there exists  in~$Q_n$ a maximum volume simplex with some vertex coinciding with a vertex of~the~cube.
 For such simplices, we have  the following theorem. 

\smallskip 
{\bf Theorem 1.} {\it Suppose $S\subset Q_n$ is a maximum volume simplex with some vertex coinciding with a vertex of the cube, 
$P$ is the corresponding interpolation projector. Then
\begin{equation}\label{th1_formula}
\|P\|_{Q_n}\leq \frac{2h_{n+1}}{h_n}+1.
\end{equation}
}

\smallskip 
{\it Proof.}
Let $x^{(1)}, \ldots, x^{(n+1)}$ be the vertices of $S$. We may assume that $x^{(n+1)}=0,$ 
i.\,e., the vertex matrix of the simplex has the form 
$${\bf S} =
\left( \begin{array}{cccc}
x_1^{(1)}&\ldots&x_n^{(1)}&1\\
\vdots&\vdots&\vdots&\vdots\\
x_1^{(n)}&\ldots&x_n^{(n)}&1\\
0&\ldots&0&1\\
\end{array}
\right).$$

Since $S$ is a maximum volume simplex in $Q_n$, we have $\vo(S)=\nu_n$ 
and~$|\Delta|=|\det({\bf S})|=n!\nu_n=h_n.$ Hence, if $x\in \ver(Q_n)$, then
$|\Delta_j(x)|\leq |\Delta|$,  $1\leq j\leq n+1.$

Fix $x\in \ver(Q_n)$ and consider determinants $\Delta_j(x)$, $1\leq j\leq n$.
From  definition of $\Delta_j(x)$ and properties of determinant, 

$$ |\Delta_j(x)|=\left| \begin{array}{cccc}
x_1^{(1)}&\ldots&x_n^{(1)}&1\\
\ldots&\ldots&\ldots&\ldots\\
x_1&\ldots&x_n&1\\
\ldots&\ldots&\ldots&\ldots\\
x_1^{(n)}&\ldots&x_n^{(n)}&1\\
0&\ldots&0&\pm1\\
\end{array}
\right | =
\left| \begin{array}{cccc}
x_1^{(1)}&\ldots&x_n^{(1)}&0\\
\ldots&\ldots&\ldots&\ldots\\
x_1&\ldots&x_n&0\\
\ldots&\ldots&\ldots&\ldots\\
x_1^{(n)}&\ldots&x_n^{(n)}&0\\
0&\ldots&0&\pm1\\
\end{array}
\right | =
$$

\medskip
$$=
\left| \begin{array}{cccc}
x_1^{(1)}&\ldots&x_n^{(1)}&0\\
\ldots&\ldots&\ldots&\ldots\\
x_1&\ldots&x_n&0\\
\ldots&\ldots&\ldots&\ldots\\
x_1^{(n)}&\ldots&x_n^{(n)}&0\\
x_1^{(j)}&\ldots&x_n^{(j)}&\pm 1\\
\end{array}
\right | =
\left| \begin{array}{cccc}
x_1&\ldots&x_n&0\\
x_1^{(1)}&\ldots&x_n^{(1)}&0\\
\ldots&\ldots&\ldots&\ldots\\
x_1^{(j)}&\ldots&x_n^{(j)}&\pm 1\\
\ldots&\ldots&\ldots&\ldots\\
x_1^{(n)}&\ldots&x_n^{(n)}&0\\
\end{array}
\right | .
$$

\medskip
\noindent
Here we mark the  $j$th row.  Consequently, there exist  numbers  $u_j=\pm 1$ such that
$$\sum_{j=1}^{n}|\Delta_j(x)|
=\left| \begin{array}{cccc}
x_1&\ldots&x_n&0\\
x_1^{(1)}&\ldots&x_n^{(1)}&u_1\\
\vdots&\vdots&\vdots&\vdots\\
x_1^{(n)}&\ldots&x_n^{(n)}&u_n
\end{array}
\right |.$$

\smallskip
\noindent
Define $n$-dimensional binary vectors $v$ and $w$ as follows. If $u_j=1$, then $v_j=1$ and~$w_j=0$;
if $u_j=-1,$ then $v_j=0$ and $w_j=1$. We have $u=v-w$, hence 

$$\sum_{j=1}^{n}|\Delta_j(x)|=
\left| \begin{array}{cccc}
x_1&\ldots&x_n&0\\
x_1^{(1)}&\ldots&x_n^{(1)}&v_1\\
\vdots&\vdots&\vdots&\vdots\\
x_1^{(n)}&\ldots&x_n^{(n)}&v_n
\end{array}
\right| \quad -
\quad \left| \begin{array}{cccc}
x_1&\ldots&x_n&0\\
x_1^{(1)}&\ldots&x_n^{(1)}&w_1\\
\vdots&\vdots&\vdots&\vdots\\
x_1^{(n)}&\ldots&x_n^{(n)}&w_n
\end{array}
\right|.$$

\smallskip
\noindent
Denote these determinants of order $n+1$ by $d_1$ and  $d_2$.  Since their entries are in~$[0,1]$, then $|d_1|, |d_2|\leq h_{n+1}$, 
$$\sum_{j=1}^{n}|\Delta_j(x)|=d_1-d_2\leq |d_1|+|d_2|\leq 2h_{n+1}.$$  Using
$|\Delta|=h_n,$ we get
$$\sum_{j=1}^{n}\frac{|\Delta_j(x)|}{|\Delta|}\leq \frac{2h_{n+1}}{h_n}.$$

It remains to include the case $j=n+1$. Recall that
$$\frac{|\Delta_{n+1}(x)|}{|\Delta|}\leq 1, \quad x\in \ver(Q_n),$$
therefore, for any vertex of the cube,
$$\sum_{j=1}^{n+1}\frac{|\Delta_j(x)|}{|\Delta|}\leq \frac{2h_{n+1}}{h_n}+1.$$
We obtain
$$
\|P\|_{Q_n}=  \max \limits_{x\in \ver(Q_n)}
\sum \limits_{j=1}^{n+1} |\lambda_j(x)|= \max \limits_{x\in \ver(Q_n)}
\sum \limits_{j=1}^{n+1} \frac{|\Delta_j(x)|}{|\Delta|}\leq \frac{2h_{n+1}}{h_n}+1.$$
\hfill$\Box$

\smallskip
{\bf Corollary 1.} 
\begin{equation}\label{h_nplus1_frac_hn_ineq}
\theta_n\leq \frac{2h_{n+1}}{h_n}+1.
\end{equation}

\smallskip
This follows immediately from  \eqref{th1_formula}.

\smallskip
{\bf Corollary 2.} {\it There exist $c, c_1>0$ does not depend on $n$ such that
\begin{equation}\label{corol2_formula}
\frac{h_{n+1}}{h_n}\geq c\sqrt{n}, \quad \frac{\nu_n}{\nu_{n+1}}\leq c_1\sqrt{n}.
\end{equation}
}

\smallskip
{\it Proof.} The author proved that $\theta_n\geq \gamma \sqrt{n}$,
 $\gamma>0$ not depends on $n$ (see \cite{nevskii_monograph}). Hence, the left-hand inequality in \eqref{corol2_formula} follows from
\eqref{h_nplus1_frac_hn_ineq}. Since $h_n=n!\nu_n$, the right-hand inequality appears from the left-hand one.
\hfill$\Box$

 \smallskip
 Let us note  more explicit forms of \eqref{h_nplus1_frac_hn_ineq}. For any $n$,
 $$\theta_n>\frac{\sqrt{n-1}}{e}.$$
 Consequently,
 $$\frac{h_{n+1}}{h_n}\geq\frac{\theta_n-1}{2}>\frac{\sqrt{n-1}}{2e}-\frac{1}{2},$$
 $$\frac{\nu_{n+1}}{\nu_n}=\frac{h_{n+1}}{(n+1)h_n}>\frac{\sqrt{n-1}}{2e(n+1)}-\frac{1}{2(n+1)}.$$

The above proof of Theorem 1 follows the approach of \cite{nevskii_mais_2007}.
 Unfortunately, later it was noticed  that the arguments in \cite{nevskii_mais_2007} perfomed for an arbitrary $n$
have a lacuna. Yet  in the case when $n+1$ is an Hadamard number Theorem 1 allows to obtain for  $\theta_n$ 
upper bounds which cannot be improved  in order of $n$. 
To show this, let~us utilize the estimates   
\begin{equation}\label{Hadamard_ineqs_arbitraty_n}
h_n\leq\frac{\left(n+1\right)^{(n+1)/2}}{2^n},   \quad n\in {\mathbb N};
\end{equation}
\begin{equation}\label{Hadamard_ineqs_even_n}
h_n\leq \frac{n^{n/2}\sqrt{2n+1}}{2^n}, \quad n \ \mbox{ is even}.
\end{equation}
The equality in  \eqref{Hadamard_ineqs_arbitraty_n} 
takes place if and only if $n+1$ is an Hadamard number. Relation \eqref{Hadamard_ineqs_arbitraty_n} was proved by Hadamard \cite{hadamard_1893}; inequality  \eqref{Hadamard_ineqs_even_n}  was obtained by Barba
\cite{Barba_1933} (see~\cite{hudelson_1996}).

\smallskip
{\bf Corollary 3.} {\it Suppose $n+1$ is an Hadamard number. If $S$ is a regular simplex with the vertices coinciding with vertices of $Q_n,$  then 
for the corresponding interpolation projector $P$
$$\|P\|_{Q_n}\leq \sqrt{2n+3}+1.$$
}

\smallskip
{\it Proof.} If $n+1$ is an Hadamard number, then the considered simplex $S$ has maximum volume among all simplices contained in  $Q_n$.
Hence,
 we can apply  to~$S$ Theorem 1:
$$\|P\|_{Q_n}\leq \frac{2h_{n+1}}{h_n}+1.$$
As it was noted above, in this case 
$$h_n=\frac{\left(n+1\right)^{(n+1)/2}}{2^n}.$$
Since $n+1$ is even, we can apply to $h_{n+1}$ inequality \eqref{Hadamard_ineqs_even_n}:
$$h_{n+1}\leq \frac{(n+1)^{(n+1)/2}\sqrt{2n+3}}{2^{n+1}}.$$
From this, it follows
$$\|P\|_{Q_n}\leq \frac{2h_{n+1}}{h_n}+1\leq \frac{2(n+1)^{(n+1)/2}\sqrt{2n+3} }
{2^{n+1}}\cdot \frac {2^n}{\left(n+1\right)^{(n+1)/2}}+1=$$
$$=\sqrt{2n+3}+1.$$
\hfill$\Box$

\smallskip
{\bf Corollary 4.} {\it If $n+1$ is an Hadamard number, then $\theta_n\leq \sqrt{2n+3}+1.$}

\smallskip
Upper bounds given in Corollaries 3 and 4  are unimprovable in order of $n$, but can be made more  precise. We will show this in next two sections.

\section{The Proof Based on Hadamard Matrices}\label{nev_s3}

Recall  (see  \cite{hudelson_1996}) that $n+1$ is an Hadamard number iff
it is possible to inscribe an $n$-dimensional regular simplex into an $n$-dimensional cube in such a way
 that the~vertices of the simplex coincide with vertices of the cube.
This can be easily shown  for~the~cube $Q_n^\prime=[-1,1]^n$. 
There is a simple correspondence between Hadamard matrices of order $n+1$ with the last column consisting of $1$'s and $n$-dimensional regular simplices whose  vertices are at the vertices of $Q_n^\prime$. Since we need this correspondence, we will dwell on it in more detail.

If $S$ is~a~regular simplex of the specified type, then its vertex matrix  ${\bf S}$  is an~Hadamard matrix of order $n+1$. Indeed,
 let  $x^{(1)}, \ldots,  x^{(n+1)}$ be the vertices of  $S$.
 Since $x^{(j)}$ coincide with vertices of the cube, the entries of $\bf S$ are  $\pm 1$ and~its last column consists of $1$'s. 
  Denote by $h^{(j)}$  the rows of  ${\bf S}$ and make sure that these vectors are pairwise orthogonal in ${\mathbb R}^{n+1}$.
  Equalities $||x^{(j)}||= \sqrt{n}$ mean that simplex $S$ is inscribed into an $n$-dimensional ball of radius $\sqrt{n}$. 
  The edge-length $d$ of~a~regular simplex and radius $R$ of the circumscribed ball satisfy the relation 
\begin{equation}\label{d_R_eq}
d=R\sqrt{2}\sqrt{\frac{n+1}{n}}
\end{equation}
 (see, e.\,g.,~\cite{nevskii_monograph}).
If $R=\sqrt{n},$ then $d=\sqrt{2(n+1)}.$  Hence,
$$2(n+1)=||x^{(j)}-x^{(k)}||^2=2||x^{(j)}||^2
- 2(x^{(j)},x^{(k)})=2n-2(x^{(j)},x^{(k)}).$$
 We see that  $(x^{(j)},x^{(k)})=-1.$ 
Since the rows of $\bf S$ in the first $n$ columns contain coordinates of vertices of the simplex and the last element of any row is $1$, then 
$$(h^{(j)},h^{(k)})=(x^{(j)},x^{(k)})+1=0.$$
Thus,  $h^{(j)}$  are pairwise orthogonal in ${\mathbb R}^{n+1}$.
This means that $\bf S$ is an Hadamard matrix of order $n+1$. We have
\begin{equation}\label{bfS_is_Hadamard}
{\bf S}^{-1}=\frac{1}{n+1}\, {\bf S}^{T}.
\end{equation}

Conversely, let ${\bf H}$ be an Hadamard matrix of order $n+1$ with the last column consisting of  $1$'s. Consider the simplex   $S$ whose vertex matrix is  $\bf H$.
In other words, the~vertices of  $S$ are given by the rows of   ${\bf H}$  (excepting the last component).  The~simplex $S$ is inscribed into the cube $Q_n^\prime,$ moreover,
$\ver(S)\subset \ver(Q_n^\prime)$. If $h^{(j)}$ are the~rows of ${\bf H}$ and $x^{(j)}$ are the vertices of  $S$, then 
$(h^{(j)},h^{(k)})=0,$  $(x^{(j)},x^{(k)})=-1,$  $\|x^{(j)}\|^2=n.$
It follows that $||x^{(j)}-x^{(k)}||^2=2(n+1)$. So,  $S$ is a regular simplex with the edge-length  $\sqrt{2(n+1)}$.

Notice that any Hadamard matrix of order  $n+1$ is equivalent to the vertex matrix of some   $n$-dimensional regular simplex inscribed in $Q_n^\prime$. This vertex matrix
can be obtained from the given Hadamard  matrix after multiplication of some rows by~$-1$.

\smallskip 
{\bf Theorem 2.} {\it  Suppose $n+1$ is an Hadamard number,  $S$  is an $n$-dimensional regular simplex having the vertices at vertices of  $Q_n^\prime$. Then for the
corresponding interpolation projector $P:C(Q_n^\prime)\to \Pi_1({\mathbb R}^n)$ 
\begin{equation}\label{th2_formula}
\|P\|_{Q_n^\prime}\leq \sqrt{n+1}.
\end{equation}
}

\smallskip 
{\it Proof.}  Let  $x^{(1)}, \ldots,  x^{(n+1)}$ be the vertices and  let $\lambda_1, \ldots,  \lambda_{n+1}$ be the basic \linebreak Lagrange polynomials of the simplex.
Under our propositions, the vertex matrix $\bf S$ is an~Hadamard matrix of order  $n+1$ with the last column consisting of $1$'s.

Let's show that for $x\in{\mathbb R}^n$ 
\begin{equation}\label{proof_th2_1}
\sum_{j=1}^{n+1} \lambda_j(x)^2=\frac{||x||^2+1}{n+1}. 
\end{equation}
Denote $y=(x_1,\ldots,x_n,1)\in {\mathbb R}^{n+1}.$ 
The coefficients of 
$\lambda_j$ form the $j$th column of~${\bf S}^{-1}$. Since  $\bf S$ is an Hadamard matrix, it satisfies~\eqref{bfS_is_Hadamard}, whence
$$\lambda_j(x)=\frac{1}{n+1}(h^{(j)},y).$$
Further,  $\displaystyle\frac{h^{(j)}}{\sqrt{n+1} }$  form an ortho\-normalized basis of ${\mathbb R}^{n+1},$ therefore,
$$y=\sum_{j=1}^m \frac{ (h^{(j)},y) }{\sqrt{n+1}}
\frac{h^{(j)}}{\sqrt{n+1}}, \quad (y, y)=\sum_{j=1}^{n+1}
\frac{(h^{(j)},y)^2} {n+1}.$$ 
From this, we have
$$||x||^2+1=(y, y)=\sum_{j=1}^{n+1}
\frac{(h^{(j)},y)^2} {n+1}=(n+1)\sum_{j=1}^{n+1} \lambda_j(x)^2.$$
Last equalities yield  \eqref{proof_th2_1}.

If $x$ is a vertex of $Q_n^\prime$, then  $||x||^2=n$ and \eqref{proof_th2_1} implies
$\sum \lambda_j(x)^2=1.$
Applying Cauchy inequality, we get for $x\in \ver(Q_n^\prime)$
$$\sum_{j=1}^{n+1}|\lambda_j(x)|\leq \left(\sum_{j=1}^{n+1} \lambda_j(x)^2\right)^{\frac{1}{2}}
\cdot\sqrt{n+1}=\sqrt{n+1}.$$
Suppose
 $P:C(Q_n^\prime)\to \Pi_1({\mathbb R}^n)$ is the interpolation projector corresponding to $S$.
Then we have
$$||P||_{Q_n^\prime}=\max_{x\in\ver(Q_n^\prime)} \sum_{j=1}^{n+1} |\lambda_j(x)|\leq \sqrt{n+1}.$$
Theorem is proved.
\hfill$\Box$

\smallskip
By similarity reasons, Theorem 2 overcomes to any $n$-dimensional cube, for instance, to the cube 
$Q_n=[0,1]^n$.

\smallskip
{\bf Corollary 5.} {\it Suppose $n+1$ is an Hadamard number. Then $\theta_n\leq \sqrt{n+1}.$}

\smallskip
The result of Corollary 5 is known
(see \cite{nevskii_mais_2003}, \cite{nevskii_monograph}), but the given proof  is more clearly related to Hadamard matrices. 

Note that $n$-dimensional regular simplices with vertices at vertices of the cube can be located differently with respect to  vertices and faces of the cube. 
This is quite evident when  the norms of the corresponding projectors are different. But this is also possible if regular simplices generate the same norms. 
In more detail, we will  describe an approach based on comparison of $\mu$-vertices of the cube with respect to~various~simplices.

A notion  {\it $\mu$-vertex of the cube $Q_n$  with respect to a simplex $S\subset Q_n$} was introduced by the author in   \cite{nevskii_mais_2009}.
 If  in place of $Q_n$ we consider an arbitraty \linebreak $n$-dimensional cube $Q$,  we obtain  equivalent results.  

Assume $1\leq \mu \leq n.$ 
We say that  $x\in\ver(Q)$ is
a
$\mu$-vertex of a cube $Q$ with respect to a simplex  $S\subset Q,$
iff for the interpolation projector  
$P:C(Q)\to\Pi_1\left({\mathbb R}^n\right)$ 
 with the~nodes at vertices of $S$ holds 
$$\sum_{j=1}^{n+1} |\lambda_j(x)|=\|P\|_{Q}$$
and among numbers $\lambda_j(x)$ there are exactly $\mu$ negatives.
This notion is closely connected with some relation between  $\|P\|_Q$ and the value
$$\xi(Q;S):=\min \{\sigma \geq 1: Q\subset \sigma S\}.$$
Here $\sigma S$ denotes the homothetic copy of   $S$ with the center of homothety at~the~center of gravity of $S$ and the ratio of homothety $\sigma$.  
We call  $\xi(Q;S)$  {\it the~absorption index of the cube  $Q$ by the simplex $S$}.
 For a projector $P:C(Q)\to\Pi_1\left({\mathbb R}^n\right)$ 
and the corresponding simplex  $S$, it is proved in  \cite{nevskii_mais_2009} that  
\begin{equation}\label{xi_normP_ineqs}
\frac{n+1}{2n}\left( \|P\|_{Q}-1\right)+1\leq \xi(Q;S) \leq 
\frac{n+1}{2}\left( \|P\|_{Q}-1\right)+1. 
\end{equation}
The right-hand equality in \eqref{xi_normP_ineqs} takes place if and only if there exists an  $1$-vertex of~$Q$ with respect to $S.$ If for some  $\mu$
there is a  $\mu$-vertex of
$Q$ with respect to $S,$ then
$$\frac{n+1}{2\mu}\left( \|P\|_{Q}-1\right)+1\leq \xi(Q;S).$$

Of course, simplices which have different sets of $\mu$-vertices with respect to~the~containing cube, are differently  located in the cube, even if have the same projectors' norms.
These simplices are non-equivalent in the following sense: one of them cannot be mapped into another by an orthogonal transform which maps the cube into~itself.
Let us give some examples concerning $Q=Q_n^\prime.$  

While obtaining estimates for minimal norms of projectors, in \cite{nev_ukh_mais_2018_25_3} various \linebreak $n$-dimensional regular simplices arising from different Hadamard matrices of~order $n+1$ were calculated.
 For~$n=15$, the order of matrices is equal to  $16$. Up~to~equivalence, there are exactly five these Hadamard matrices. They correspond to five simplices described in Table \ref{tab:nev_Reg_simplices_for_n_15}. 
 The results of  calculations were kindly delivered to~the author by A.\,Yu. Ukhalov.

\begin{table}[!htbp]
\begin{center}
\caption{\label{tab:nev_Reg_simplices_for_n_15} Regular simplices for $n=15$}
\medskip
\bgroup
$
\def\arraystretch{1.5}
\begin{array}{|c|c|c|c|}
  \hline
 S & {\|P\|}_{Q_{15}^\prime} & \mbox{\bf Values } \ \mu& \mbox{\bf Number of} \ \mu\mbox{\bf-vertices}  \\ 
  \hline
 S_1 & 4 & 6& m_6=448 \\
  S_2 & 4& 6& m_6=192 \\ S_3 & 4& 6& m_6=64 \\ S_4 & \frac{7}{2} & 4, 5, 6, 8& m_4=896, m_5=1344, m_6=5376, m_8=1344 \\
   S_5 & \frac{7}{2} & 4, 5, 6, 8& m_4=896, m_5=1344, m_6=5376, m_8=1344 \\ 
  \hline
 \end{array}
$
\egroup
\end{center}
\end{table}

By $m_\mu$ we denote a number of $\mu$-vertices of the cube  $Q_n^\prime$ with respect to every simplex. For the rest 
$1\leq\mu\leq 15$, 
excepting the given in Table \ref{tab:nev_Reg_simplices_for_n_15},
values  $m_\mu$ are~zero.
Every simplex $S_1$, $S_2$, and $S_3$ generates the same projector norm and has only \linebreak $6$-vertices. But the numbers of $6$-vertices for them are different,
hence, these simplices are pairwise non-equivalent. Each of them also is non-equivalent  as to $S_4$, so as to $S_5.$
The latter simplices generates equal norms and have the same sets \linebreak of $\mu$-vertices. Also we have $\theta_{15}\leq \frac{7}{2}=3.5$.
This is more exact than the estimate  $\theta_{15}\leq 4$   of Corollary 5 for $n=15$.

Another example is related to $n=23.$ In \cite{kudr_ozerova_ukh_2017}, the results for 60 regular simplices inscribed into $Q_{23}^\prime$ were mentioned. 
The simplices were built from the available 60~pairwise non-equivalent Hadamard matrices of order  $24$. For all simplices, excepting ones with the numbers 16, 53, 59, and 60,
the projector norm  equals $\frac{14}{3}=4.6666\ldots$, while for~each of these four simplices the projector norm is $\frac{9}{2}=4.5$. In~particular, we have  $\theta_{23}\leq 4.5$
(what is noted also  in  \cite{nev_ukh_mais_2018_25_3}). This inequality is more exact than the estimate  $\theta_{23}\leq \sqrt{24}=4.8989\ldots$  of~Corollary 5 for $n=23$.  Each  of the four exceptional simplices is not equivalent to any of the 56 others.

Despite the possible differs, for all regular simplices with verices at vertices \linebreak of the~cube, inequality \eqref{th2_formula} holds. The inscribed regular simplices satisfying $\|P\|_{Q_n^\prime}=\sqrt{n+1}$ exist at least for $n=1$, $n=3$, and $n=15.$ The problem of full description of dimensions $n$ with such property is still open.



\section{Connection with Interpolation on a  Ball}\label{nev_s4}

A regular simplex inscribed into an $n$-dimensional ball has maximal volume among all simplices contained in this ball.   There are no another simplices \linebreak having this
property  (see \cite{fejes_tot_1964}, \cite{slepian_1969}, and
\cite{vandev_1992}). In the case when  $n+1$ is an Hadamard number, the similar property with respect to an $n$-dimensional cube holds for a~re\-gular simplex inscribed into the cube.  
Furthermore, the upper bound for  projectors' norm  corresponding to  these simplices   both on a cube and on a ball  is equal
 to~$\sqrt{n+1}$. 

In this section, we give another proof of Theorem 2. It seems to be shorter than previous one from Section 3, but  is based on some relations  which were obtained in \cite{nev_ukh_mais_2019_2}  in a rather complicated way. This approach is related with interpolation on~a~Euclidean ball.

At first, suppose $S$ is a regular simplex inscribed into an $n$-dimensional ball  \linebreak
$B=B(x^{(0)};R)$ and 
$P:C(B)\to \Pi_1({\mathbb R}^n)$  is the corresponding interpolation projector. Obviously, $\|P\|_B$ does not depend neither on the center $x^{(0)}$ and the radius $R$
of~the~ball, nor on the chosen  regu\-lar simplex inscribed into this ball. In other words,
$\|P\|_B$  depends only on dimension $n$. 
For $0\leq t\leq n+1$, consider the function
\begin{equation}\label{psi_function_modulus}
\psi(t):=\frac{2\sqrt{n}}{n+1}\Bigl(t(n+1-t)\Bigr)^{1/2}+
\left|1-\frac{2t}{n+1}\right|.
\end{equation}
Let us denote
$\displaystyle a:=\left\lfloor\frac{n+1}{2}-\frac{\sqrt{n+1}}{2}\right\rfloor.$
It was proved in~\cite{nev_ukh_mais_2019_2} that
$$\|P\|_B=\max\{\psi(a),\psi(a+1)\}.$$
This equality yields
\begin{equation}\label{norm_P_reg_ineqs}
\sqrt{n}\leq \|P\|_B\leq \sqrt{n+1}.
\end{equation}
Moreover, $\|P\|_B=\sqrt{n}$ only when $n=1$, and 
$\|P\|_B=\sqrt{n+1}$ if and only if
$\sqrt{n+1}$ is an integer.

 Theorem 2 follows immediately from the right-hand inequality \eqref{norm_P_reg_ineqs}. 
 Assume $n+1$ is an Hadamard number.
Consider an interpolation projector $P$ with the~nodes at~those vertices of the cube $Q_n^\prime$ that form a regular simplex $S$. Since  $Q_n^\prime$ is inscribed into the unit
ball $B_n$,   simplex $S$ is also inscribed into $B_n$. It remains to apply  formula~\eqref{norm_of_P_lambda_j} for the projector's norm both on the cube and on the ball and
also  inequality \eqref{norm_P_reg_ineqs}  in~the case $B=B_n$:
$$\|P\|_{Q_n^\prime}=\max_{x\in Q_n^\prime}\sum_{j=1}^{n+1}
|\lambda_j(x)|\leq \max_{x\in B_n}\sum_{j=1}^{n+1}
|\lambda_j(x)|= \|P\|_{B_n}\leq \sqrt{n+1}.$$
Inequality \eqref{th2_formula} is proved.

The upper bound $\sqrt{n+1}$ of the projector's norm is delivered in interpolation on~the ball or on the cube in different ways.
For dimensions $n=m^2-1$ , and only in these cases,  
the equality $\|P\|_{B_n}=\sqrt{n+1}$ holds for any
regular simplex inscribed into
$B_n$. If $n+1$ is an Hadamard number, then the equality
 $\|P\|_{Q_n^\prime}=\sqrt{n+1}$ may hold as for all regular simplices having vertices at vertices of the cube $(n=1, n=3)$, as for some of them $(n=15),$
 or  may not be executed at all.

The connection with constructions on a ball can be also seen in the proof  given in Section 3 (see \eqref{d_R_eq}).



\begin{thebibliography}{99}
\bibitem{Barba_1933}
{\sc G. Barba.} Intorno al teorema di Hadamard sui determinanti a valore massimo, {\it Glornale Mat. Battaglini (3),} {\bf 71}\,(1933), 70--86.

\bibitem{fejes_tot_1964}
 {\sc L. Fejes T\'{o}t.}
{\it Regular Figures},
 New York: Macmillan/Pergamon, 1964.

\bibitem{hadamard_1893}
{\sc J.  Hadamard.}  R\'esolution d'une question relative aux d\'eterminants,
{\it Bull. Sciences Math. (2),} {\bf 17}\,(1893), 240--246.

\bibitem{hall_1970} 
{\sc M. Hall, Jr.} {\it Combinatorial Theory,} Blaisdall Publishing Company, Waltham (Massachusetts)--Toronto--London, 1967.

\bibitem{horadam_2007}
{\sc K.\,J. Horadam.} {\it Hadamard Matrices and Their Applications,} Princeton University Press, 2007.


\bibitem{hudelson_1996} 
 {\sc M. Hudelson, V. Klee,  and D. Larman.} Largest $j$-simplices
in $d$-cubes: some relatives of the Hadamard maximum determinant
problem, {\it Linear Algebra Appl.,} {\bf 241--243}\,(1996), 519--598.

\bibitem{kudr_ozerova_ukh_2017} 
{\sc I.\,S. Kudryavcev,  E.\,A. Ozerova, and  A.\,Yu. Ukhalov.}  New estimates for the norms of minimal projectors,
  In:\,{\it Sovremennye Problemy Matematiki i~Informatiki. Vypusk 17} (Modern Problems in Mathematics and Informatics. Issue 17),
    P.\,G. Demidov Yaroslavl  State University, 
  Yaroslavl, 2017, 74--81
 (in~Russian). 

\bibitem{manjhi_2022}
{\sc P.\,K. Manjhi and M.\,K. Rama.} Some new examples of circulant partial Hadamard matrices of type $4 - H(k\times n)$,
{\it Advances and Applications in Mathe\-matical Sciences,} {\bf 21}:5\,(2022), 2559--2564.

\bibitem{nevskii_mais_2003}
{\sc M.\,V. Nevskii.} Estimates for minimal norm of a projector in linear  inter\-polation over vertices of an $n$-dimensional cube, {\it  Modeling and Analysis \linebreak of Information Systems,}
{\bf10}:1\,(2003), 9--19.

\bibitem{nevskii_mais_2007}
{\sc M.\,V. Nevskii.} Minimal projectors and largest simplices, {\it  Modeling and Ana\-lysis of~Information Systems,}
{\bf14}:1\,(2007), 3--10.

\bibitem{nevskii_mais_2009}
{\sc M.\,V. Nevskii.} On a certain relation for the minimal norm of an interpolation projector, {\it  Modeling and Ana\-lysis of~Information Systems,}
{\bf16}:1\,(2009), 24--43.



\bibitem{nevskii_monograph}
{\sc M.\,V. Nevskii.}
{\it Geometricheskie Ocenki v Polinomialnoi Interpolyacii} 
(Geometric Estimates in Polynomial
Interpolation), Yaroslavl': Yarosl. Gos. Univ., 2012 (in~Russian).


\bibitem{nev_ukh_mais_2018_25_3}
{\sc M.\,V. Nevskii and  A.\,Yu. Ukhalov.}
On optimal interpolation by linear functions on an $n$-dimensional cube, {\it  Modeling and Ana\-lysis of~Information Systems,}
{\bf 25}:3\,(2018), 291--311

\bibitem{nev_ukh_mais_2019_2}
{\sc M.\,V. Nevskii and A.\,Yu. Ukhalov.} 
 Linear interpolation on a Euclidean ball in~${\mathbb R}^n$,
{\it  Modeling and Ana\-lysis of~Information Systems,} {\bf 26}:2\,(2019), 279--296.











 
\bibitem{slepian_1969}
{\sc D. Slepian.}
The content of some extreme simplices,
{\it Pacific J.~Math.,} {\bf 31}\,(1969), 795--808.

\bibitem{vandev_1992}
{\sc D.  Vandev.}
A minimal volume ellipsoid around a simplex,
{\it C. R. Acad. Bulg.~Sci.}, {\bf 45}:6\,(1992), 37--40.



\end{thebibliography}
\end{document}